\documentclass[11pt]{article}
\usepackage{fullpage}
\usepackage{amsfonts, amsmath, amssymb,latexsym,epsfig}
\usepackage{amsthm}
\usepackage{tikz}
\usetikzlibrary{graphs}
\usetikzlibrary{graphs.standard}
\usepackage{relsize}
\usepackage{verbatim}
\usepackage{enumerate}
\usepackage[skip=10pt plus1pt, indent=20pt]{parskip}

\newtheorem{thm}{Theorem}[section]
\newtheorem{cor}[thm]{Corollary}
\newtheorem{lem}[thm]{Lemma}

\newtheorem{prop}[thm]{Proposition}

\newcommand{\F}{\mathbb F}

\providecommand{\keywords}[1]{\textbf{\textit{Keywords}} #1}

\title{Achromatic colorings of polarity graphs}
\author{Vladislav Taranchuk\thanks{Department of Mathematical Sciences, University of Delaware. \texttt{vladtar@udel.edu}} \and Craig Timmons\thanks{Department of Mathematics and Statistics, California State University Sacramento, 6000 J Street, Sacramento, CA 95819.   
		\texttt{craig.timmons@csus.edu}}}

\begin{document}
	
	\maketitle

	\begin{abstract}
		A complete partition of a graph $G$ is a partition of the vertex set such that there is at least one edge between any two parts.  The largest $r$ such that $G$ has a complete partition into $r$ parts, each of which is an independent set, is the achromatic number of $G$.  We determine the achromatic number of polarity graphs
		of biaffine planes coming from generalized polygons.
		Our colorings of a family of unitary polarity graphs 
		are used to solve a problem of Axenovich and Martin on complete partitions of $C_4$-free graphs.  Furthermore, these colorings prove that there are sequences of graphs which are optimally complete and have unbounded degree, a problem that had been studied for the sequence of hypercubes  independently by Roichman, and Ahlswede, Bezrukov, Blokhuis, Metsch, and Moorhouse.      
	\end{abstract}
	
	\keywords{polarity graph, achromatic coloring, complete partitions}
	
	\section{Introduction}
	
	A \emph{complete partition} of a graph $G$ is a partition of $V(G)$ into parts $V_1 , \dots , V_r$ such that there is at least one edge between any two parts.  Equivalently, a \emph{pseudo-achromatic coloring} of $G$ is a vertex coloring such that there is an edge between every pair of distinct color classes.    
	The largest $r$ such that $G$ has a pseudo-achromatic coloring with $r$ colors is the 
	\emph{pseudo-achromatic number of $G$}, denoted $\psi (G)$.  Different notation has been used for $\psi (G)$, in part because it appears in different areas.  For example, $\tau (G)$ was used by Kostochka \cite{Kostochka} and $\textup{ccl}' (G)$ was used by Bollob\'{a}s, Erd\H{o}s, and Catlin \cite{Bollobas} in their work 
	on Hadwiger's conjecture and random graphs. 
	We are following the notation and terminology of
	Thomason \cite{Thomason}, which is the same notation used by Yegnanarayanan \cite{Yeg1, Yeg2, Yeg3}, who has written several papers on pseudo-achromatic colorings.  
	
	A pseudo-achromatic coloring where the color classes are independent sets is an \emph{achromatic coloring}.  The \emph{achromatic number} of $G$, denoted 
	$\chi_a(G)$, is the largest $r$ such that $G$ has an 
	achromatic coloring with $r$ colors.  
	For any graph $G$, $\chi (G) \leq \chi_a (G) \leq \psi (G)$
	where $\chi (G)$ is the chromatic number of $G$.  It is interesting
	to study the differences in these graph parameters.  For example, Edwards \cite{EdwardsTrees} gave an infinite family of trees $\{ T_r \}_{ r = 1 }^{ \infty}$ such that $\psi (T_r) = \chi_a (T_r) + 1$.  This disproved the conjecture of Hedetniemi (see \cite{EdwardsTrees}) that asserted $\psi (T)  = \chi_a (T)$ for every tree $T$.  The conjecture was later recast as $\chi_a (T) \leq \psi (T) \leq \chi_a (T) + 1$ 
	in \cite{Yeg3} Conjecture 3, but we are not aware of additional developments.  
	
	In a complete partition of $G$ into $r$ parts, there is at least one edge between each pair of  parts, implying $\binom{r}{2} \leq e(G)$.  This implies $\chi_a (G) \leq \sqrt{ 2 e(G) +1/4} + 1/2$ and when equality holds, there is exactly one edge between every pair of parts and no edge within any part.  Roichman \cite{Roichman} called such a coloring of a graph \emph{optimally complete}.  Observe that $\chi_a(G) = \psi (G)$ if $G$ is optimally complete.   Because an optimally complete coloring requires that the number of edges of $G$ be of the form $\binom{r}{2}$ for some $r \in \mathbb{N}$, most graphs do not have an optimally complete coloring.  A sequence of graphs 
	$\mathcal{G} = (  G_1, G_2, \dots )$ with $  e(G_i)  \rightarrow \infty$ is called \emph{almost optimally complete} if 
	\[
	\lim_{i \rightarrow \infty} \frac{ \chi_a (G_i ) }{ \sqrt{ 2 e (G_i) } } = 1.
	\]
	The sequence of cycles $(C_3, C_4, \dots )$ and paths 
	$(P_1, P_2, P_3, \dots )$ are almost optimally complete \cite{Yeg1, Yeg2}.  More generally, Cairnie and Edwards \cite{Cairnie1} proved that any sequence of graphs whose degrees are bounded is almost optimally complete.  In 2000, Roichman \cite{Roichman} studied $\chi_a (Q_d)$ where $Q_d$ is the $d$-dimensional hypercube, and proved that \[
	\left( \frac{1}{ 4 \sqrt{2} } + o(1) \right)  \sqrt{ 2 e (Q_d) }  \leq 
	\chi_a (Q_d) \leq \sqrt{ 2 e (Q_d) +1/4}  + 1/2.
	\]
	Here the upper bound follows from $\binom{ \psi (Q_d) }{2} \leq e(Q_d)$.  
	Less than a decade earlier, Ahlswede, Bezrukov, Blokhuis, Metsch, and Moorhouse \cite{Ahlswede} proved bounds on $\psi ( Q_d)$.  Their work is phrased in terms of vertex partitions of the hypercube into sets of pairwise distance 1.  For a sequence of graphs $\mathcal{G}= (G_1,G_2, \dots)$,
	let 
	\[
	\underline{ \psi } ( \mathcal{G} ) 
	=  \liminf_{i \rightarrow \infty} 
	\frac{\psi (G_i ) } { \sqrt{ 2 e(G_i) }} 
	~~ \mbox{and} ~~
	\underline{ \chi_a } ( \mathcal{G} ) 
	=  \liminf_{i \rightarrow \infty} 
	\frac{\chi_a (G_i ) } { \sqrt{ 2 e(G_i) }}.
	\]
	Ahlswede et al.\ noted that sequences of random graphs can be used to give examples for which these limits are 0 (see \cite{Bollobas, Kostochka}).  They asked for explicit families $\mathcal{G}$ for which
	$\underline{ \psi}( \mathcal{G} )  = 0$ and suggested that this could be true for Payley graphs.  For a prime power $q$ with $q \equiv 1 ( \textup{mod}~4)$, let $P_q$ be the Payley graph on $\mathbb{F}_q$.    
	Maistrelli and Penman \cite{Maistrelli} proved that
	$\chi_a (P_q) = \frac{q+1}{2}$.  They stated that their argument is based on Thomason's \cite{Thomason}
	proof that the Hadwiger number of $P_q$ is $\frac{q+1}{2}$.  Recall that the Hadwiger number of a graph $G$ is the largest $r$ such that $K_r$ is a minor of $G$, and this number is bounded above by $\psi (G)$.  Therefore, for the sequence of Payley graphs $\mathcal{P}_q = (P_q)$ and hypercubes $\mathcal{Q}_d = (Q_d)$, the results of \cite{Ahlswede, Maistrelli, Roichman, Thomason} imply
	\[
	\underline{ \psi } ( \mathcal{P}_q ) = \frac{\sqrt{2}}{2} ~~~
	\mbox{and} ~~~
	\frac{1}{4 \sqrt{2}} \leq
	\underline{ \chi_a } ( \mathcal{Q}_d ) \
	~~~
	\mbox{and}~~~
	\frac{\sqrt{2}}{ 2} \leq \underline{\psi } ( \mathcal{Q}_d ).
	\]
	
	As indicated by the history, various paths have led research groups to complete partitions.  Initially our work was motivated by a problem from extremal graph theory, which we will introduce in a moment.  For now, let us state one of our first results which gives a family of optimally complete graphs.  For each prime power $q$, let $G_{q^2}^{ \pi }$ be the graph with vertex set $\mathbb{F}_{q^2} \times \mathbb{F}_{q^2}$ where $(x,y)$ is adjacent to $(z,w)$ if 
	$x z^q = y + w^q$.  The graph $G_{q^2}^{ \pi }$ is a unitary polarity graph of the biaffine part of the classical projective plane, and we can determine its achromatic number exactly.  
	
	\begin{thm}\label{thm:achromatic biaffine}
		For each prime power $q$, there is a polarity graph $G_{q^2}^{ \pi }$ of the biaffine plane in $PG(2,q^2)$ such that 
		\[
		\chi_a ( G_{q^2}^{ \pi} ) = q^3.
		\]    
	\end{thm}
	
	In Section \ref{sec:C4} it will be shown that 
	$e(G_{q^2}^{ \pi } ) = \binom{q^3}{3}$.  By Theorem \ref{thm:achromatic biaffine}, this means that that $G_{q^2}^{ \pi }$ is optimally complete.
	Each infinite sequence of prime powers $q$ gives a sequence
	$\mathcal{G}_{q^2}^{ \pi}$ for which 
	$\underline{ \psi } ( \mathcal{G}_{q^2}^{ \pi}  ) = \underline{ \chi_a } ( \mathcal{G}_{q^2}^{ \pi}  ) = 1$.  
	
	The coloring used to prove Theorem \ref{thm:achromatic biaffine} was discovered in an attempt to solve a problem in extremal graph theory. 
	Given two graphs $G$ and $H$, the graph $G$ is $H$-free if $G$ does not contain $H$ as a subgraph.  The \emph{Tur\'{a}n number of $H$}, denoted $\textup{ex}(n,H)$ is the maximum number of edges in an $n$-vertex $H$-free graph $G$.  
	If $G$ is an $H$-free $n$-vertex graph, 
	then $\psi (G) \leq \sqrt{ 2  \textup{ex}(n,H)+1/4} + 1/2$.  This relates extremal numbers to complete partitions in a simple, but useful way, and leads to interesting partition problems for $H$-free graphs.  One is to fix the number of parts $r$, and then find an $H$-free graph $G$ with $\psi (G) \geq r$ where all color classes are small (otherwise one could connect parts with vertex disjoint edges).  If $G$ has a complete partition into $r$ parts each of size at most $k$, then $G$ is called an \emph{$(r,k)$-graph}.  Following Barbanera \cite{Barbanera}, Axenovich and Martin \cite{Axenovich}, for family of graphs $\mathcal{H}$ 
	let 
	\[
	f(r,\mathcal{H}) = \min \{ k \in \mathbb{N} : 
	\mbox{there is an $\mathcal{H}$-free $(r,k)$-graph} \}.
	\]
	Observe that if there is an $\mathcal{H}$-free $(r,k)$-graph $G$, then  
	\begin{equation}\label{ineq:trivial turan ub}
		\binom{r}{2} \leq e(G) \leq \textup{ex}(rk,\mathcal{H}).
	\end{equation}
	This inequality shows that if an upper bound on the Tur\'{a}n number 
	$\textup{ex}(N,\mathcal{H})$ is known, then (\ref{ineq:trivial turan ub}) 
	yields a lower bound on $f(r,\mathcal{H})$.  In turn, if an $\mathcal{H}$-free graph is known to give a good lower bound on $\textup{ex}(N,\mathcal{H})$, then one idea is to find a  pseudo-achromatic coloring of this $\mathcal{H}$-free graph to get an upper bound on $f(r,\mathcal{H})$.  Axenovich and Martin \cite{Axenovich} successfully used this approach to prove that for all $d \geq 1$, 
	\[
	\frac{r^{1/3}}{d^{1/3}} - o(r^{1/3} ) \leq f( r , K_{2,d+1 } ) \leq 
	2 r^{1/3} + o(r^{1/3}).
	\]
	The constant in the upper bound was improved in \cite{Byrne} for 
	$d > 1$ and this gave the correct dependence on $d$.  
	Currently the best known bounds for $K_{2,2}=C_4$ that we are aware of are
	$ r^{1/3} + o (r^{1/3})  \leq f(r,C_4) \leq 2 r^{1/3} + o (r^{1/3})$.  
	As a corollary to Theorem \ref{thm:achromatic biaffine}, we improve the upper bound by a factor of 2 giving an asymptotic formula for $f(r,C_4)$.     
	
	\begin{cor}\label{cor:C4 partition theorem}
		We have $f(r,C_4) = r^{1/3} + o(r^{1/3})$.
	\end{cor}

	Using polarity graphs related to the regular generalized quadrangle $Q(4, q)$ and split Cayley hexagon $H(q)$, we are able to prove results using similar partitions utilized in the proof Theorem \ref{thm:achromatic biaffine}.

	\begin{thm}\label{thm:c6 and c10}
		(a) Let $e \geq 1$ be an integer and $q = 2^{ 2e + 1}$.  There is a polarity graph $GQ_q^{ \pi}$ of the biaffine part of the generalized quadrangle $Q(4, q)$ of order $q$ with  
		\[
		\chi_a ( GQ_q^{ \pi } ) = q^2 .  
		\]
		
		\noindent
		(b) Let $e \geq 1$ be an integer and  
		$q = 3^{ 2e + 1}$.  There is a polarity graph $GH_q^{ \pi }$ of the  biaffine part of the split Cayley generalized hexagon of order $q$ with 
		\[
		\chi_a ( GH_q^{ \pi } )  = q^3 .  
		\]
	\end{thm}
	
	Computations in Section \ref{sec:C6 and C10} will show that $ e( GH_q^{ \pi } ) = \binom{ \chi_a ( GH_q^{ \pi } ) } {2}$ and 
	$ e( GQ_q^{ \pi } ) = \binom{ \chi_a ( GQ_q^{ \pi } ) } {2}$.  Also, like Theorem \ref{thm:achromatic biaffine} the colorings used to prove Theorem \ref{thm:c6 and c10} imply  results on $f(r, \mathcal{H})$ that 
	improve bounds from \cite{Byrne} by factor of 2 in certain cases.

	\begin{cor}\label{cor:C6 and C10}
		For infinitely many $r$ and $\rho$, 
		\begin{center}
			$f(r,\{ C_4 , C_6\} ) = ( 1 + o(1))r^{1/2}$ ~~ 
			and ~~
			$f( \rho , \{ C_4 , C_6 , C_8 , C_{10} \} ) = (1 + o(1)) \rho^{2/3} $.
		\end{center}
	\end{cor}
	
	The significance of Corollaries \ref{cor:C4 partition theorem} and \ref{cor:C6 and C10} is that they show the bound coming from (\ref{ineq:trivial turan ub}) gives the correct asymptotic 
	when $\mathcal{H} \in \{ \{ C_4 \}, \{C_4, C_6\} , \{C_4,C_6,C_8, C_{10 } \} \}$.  The order of magnitude of $f (r, \mathcal{H} )$ for these families, as well as other involving complete bipartite graphs, was determined in \cite{Byrne}.  For more on the function $f(r, \mathcal{H})$, see \cite{Axenovich} and \cite{Byrne}.  
	
	We prove Theorem \ref{thm:achromatic biaffine} and Corollary \ref{cor:C4 partition theorem} in Section \ref{sec:C4}.
	Theorem \ref{thm:c6 and c10} and Corollary \ref{cor:C6 and C10} will be proved in Section \ref{sec:C6 and C10}, while Section \ref{section:conclusion} contains some concluding remarks.

	\subsection{Notation and Terminology}
	
	The objects which inspired the work in this paper are primarily geometric objects. However, we work in the setting of algebraically defined graphs as this allows us to provide a clear and concise description of the complete partitions we define.  
	Let $\mathcal{P}$ and $\mathcal{L}$ be two copies of 
	$\mathbb{F}_q^{m}$, the $m$-dimensional vector space over the finite field $\mathbb{F}_q$ of order $q$. Call the set $\mathcal{P}$ \emph{points} and $\mathcal{L}$ \emph{lines}, with the distinction in notation by $(a) \in \mathcal{P}$ and $[a] \in \mathcal{L}$.
	For each $2 \leq i \leq m$, let $f_i : \mathbb{F}_q^{2i - 2} \rightarrow \mathbb{F}_q$ be a function.
	Define $\Gamma_q = \Gamma_q(f_2, f_3, \dots, f_m)$ to be the bipartite graph with parts $\mathcal{P}$ and $\mathcal{L}$ and where $(p) = (p_1, \dots, p_{m})$ is adjacent to $[\ell] = [\ell_1, \dots, \ell_{m}]$ if and only if  
	\begin{align*}
		\ell_2 + p_2 &= f_2(\ell_1, p_1), \\
		\ell_3 + p_3 &= f_3(\ell_1, p_1, \ell_2, p_2), \\
		\vdots \\
		\ell_{m} + p_{m} &= f_m(\ell_1, p_1, \ell_2 , p_2, \dots, \ell_{m-1}, p_{m-1}).
	\end{align*}
		
	A \emph{polarity} of an algebraically defined graph $\Gamma_q$ is an automorphism $\pi$ of $\Gamma_q$ such that $\pi^2$ is the identity map, 
	$\pi(\mathcal{P}) = \mathcal{L}$, and $\pi(\mathcal{L}) = \mathcal{P}$. If $\Gamma_q$ has a polarity $\pi$, then define a \textit{polarity graph} of $\Gamma_q$, denoted $\Gamma_q^\pi$, as follows. The vertex set of $\Gamma_q^\pi$ is $\mathcal{P}$.  Two vertices $p$ and $r$ are adjacent in $\Gamma_q^\pi$ if and only if the point $r$ is adjacent to the line $\pi(p)$ in $\Gamma_q$. We call a vertex of $\Gamma_q^{ \pi} $ an \textit{absolute point} if $p$ is adjacent to $\pi(p)$ in $\Gamma_q$.   
	
	In \cite{LUW}, Lazebnik, Ustimenko, and Woldar describe an important relationship between the incidence graph of a geometry and a polarity graph of the same geometry (if one exists). Since every algebraically defined graph is in fact the incidence graph of some geometry, the properties translate directly to algebraically defined graphs and to their respective polarity graphs. The most important property from \cite{LUW} that yields Corollary \ref{cor:C4 partition theorem} and Corollary \ref{cor:C6 and C10} is that if an algebraically defined graph $\Gamma_q$ does not contain an even cycle $C_{2k}$, then neither does any polarity graph of $\Gamma_q$. We include the theorem here, worded in the language of algebraically defined graphs.  
	
	\begin{thm}[{\cite[Theorem 1]{LUW}}]
		Let $\pi$ be a polarity of the algebraically defined graph $\Gamma_q$ and let $\Gamma_q^\pi$ be the corresponding polarity graph.
		\begin{enumerate}
			\item $\deg_{\Gamma_q^\pi}(p) = \deg_{\Gamma_q}(p) - 1$ if $p$ is an absolute point of $\Gamma_q^{ \pi }$. Otherwise, $\deg_{\Gamma_q^\pi}(p) = \deg_{\Gamma_q}(p)$.
			
			\item $|V(\Gamma_q^\pi)| = \frac{1}{2}|V(\Gamma_q)|$ and $|E(\Gamma_q^\pi)| = |E(\Gamma_q)| - N_\pi$, where $N_\pi$ is the number of absolute points of $\Gamma_q^{ \pi }$.
			\item If $\Gamma_q^\pi$ contains a $C_{2k+1}$, then $\Gamma_q$ contains a $C_{4k+2}$.
			\item If $\Gamma_q^\pi$ contains a $C_{2k}$, then $\Gamma_q$ contains two vertex disjoint $C_{2k}$'s. Consequently, if $\Gamma_q$ is $C_{2k}$-free, then so is $\Gamma_q^\pi$.
			\item The girth of $\Gamma_q^\pi$ is at least half the girth of $\Gamma_q$.
		\end{enumerate}
	\end{thm}

The algebraically defined graphs we use are based on generalized $n$-gons.  Let $n \geq 2$.  A \emph{generalized $n$-gon} $G$ is point-line incidence geometry whose corresponding incidence graph
$\Gamma_G$ 
has diameter $n$ and girth $2n$.  The $n$-gon $G$ is said to have order $(s,t)$ if in $\Gamma_G$, the vertices in one part have degree $s+1$ and the vertices in the other part have degree $t+1$.   
We will restrict our discussion to the case $s = t$ as these are the only geometries we deal with in this paper.

If $G$ is a generalized $n$-gon of order $(q,q)$, then 
$\Gamma_G$ has $2(q^n + q^{n-1} + \cdots +q + 1)$ vertices and $(q+1)(q^n + q^{n-1} + \cdots + q + 1)$ edges. Such geometries it turns out are not so easy to come by, and are only known to exist for $n = 2, 3, 4, 6$. Tits [CITE] proved that in fact generalized $n$-gons of order $(s, t)$ can only exist when $n = 2, 3, 4, 6$ and $8$.


	It is known that the incidence graphs of the biaffine parts of projective planes, the classical generalized quadrangle and the classical generalized hexagon can be represented as algebraically defined graphs. The exact equations we will be using for these graphs can be found in Section 4.2 of the survey of Lazebnik, Sun and Wang \cite{LSW}.  Every algebraically defined graph can be seen as the incidence graph of a geometry. 

        A generalized $n$-gon $G$ of order $(s, t)$ is an incidence geometry of points and lines whose incidence graph $\Gamma_G$ is a biaprtite biregular graph with degrees $s+1$ and $t+1$ having with girth $2n$ and diameter $n$. 
        
        When $n = 3$, the geometry is called a projective plane, when $n = 4$ the geometry is called a generalized quadrangle and when $n = 6$, the geometry is called a generalized hexagon. The most prominent examples of infinte families of each case are as follows:
        \begin{itemize}
            \item For $n = 3$, we have the Desarguesian projective plane [CITE].
            \item For $n = 4$, we have the generalized quadrangle $Q(4, q)$ (and its dual) [CITE].
            \item For $n = 6$, we have the split Cayley hexagon $H(q)$ (and its dual) [CITE].
        \end{itemize}

        An important object related to generalized $n$-gons is their corresponding \textit{biaffine part}. From a graph theoretic view, the incidence graph of the biaffine part of a generalized $n$-gon can be understood as an induced subgraph of the incidence graph of the generalized $n$-gon. Let $GP$ denote a generalized $n$-gon and $\Gamma_{GP}$ its point-line incidence graph. The biaffine part can obtained in the following manner:
        \begin{enumerate}
            \item Choose any edge in $(p, \ell) \in E(\Gamma_{GP})$, this corresponds to a flag $(p, \ell)$ in $GP$.
            \item Define the set $S$ to be the set of all vertices at distance $n -2$ from either $p$ or $\ell$ in $\Gamma_{GP}$. Since the girth of $\Gamma_{GP}$ is $2n$, we can see that $|S| = 2(q^{n-1} + q^{n-2} + \cdots + q +1)$ and that $\Gamma_{GP}$ induced on $S$ is a tree centered at the edge $(p, \ell)$.
            \item  Finally,  $\Gamma_{GP}$ induced on $V(\Gamma_{GP})\setminus S$ is called a biaffine part of $\Gamma_{GP}$. A biaffine part of a regular generalized $n$-gon has is $q$-regular, has $2q^n$ vertices, and  $q^{n+1}$ edges.
        \end{enumerate}

        The definition above implies that one generalized $n$-gon may very well have many different and non-isomorphic biaffine parts. However, in the case of the Desarguesian projective plane, the generalized quadrangle $Q(4, q)$, and split Cayley hexagon, the geometries are known to be flag-transitive (implying their corresponding incidence graphs are edge-transitive) and thus any two biaffine parts coming from one of these geometries will be isomorphic. Thus we refer to \textit{the} biaffine part of these geometries. The biaffine part of any of these mentioned geometries can be represented as an algeberaically defined graph.

	
	\section{Proof of Theorem \ref{thm:achromatic biaffine}}\label{sec:C4} 
	
	Let $q$ be a power of a prime. Denote by $G_{q^2}$ the algebraically defined graph over $\F_{q^2}^2$ whose adjacency relation is defined by the single equation
	$$
	p_2 + \ell_2 = p_1\ell_1.
	$$
	The graph $G_{q^2}$ is isomorphic to the 
	incidence graph of the biaffine part of the Desarguesian projective plane of order $q^2$.  
	Let $\pi: V(G_{q^2}) \rightarrow V(G_{q^2})$ be the function defined by 
	\begin{center}
		$\pi((p_1, p_2)) = [p_1^q, p_2^q]$ ~~~ and ~~~
		$\pi([\ell_1, \ell_2]) = (\ell_1^q, \ell_2^q)$.
	\end{center}
	Using the identities $x^{q^2} = x$ and $(x+y)^q = x^q + y^q$  for all $x,y \in \mathbb{F}_{q^2}$, it can be checked directly that 
	$\pi$ is a polarity.  Alternatively, $\pi$ is equivalent to a restriction of the classical unitary polarity 
	
	\begin{center}
		$(p_1,p_2,p_3) 
		\rightarrow [p_1^q , p_2^q , p_3^q]$,~~~~~ $[ \ell_1 , \ell_2 , \ell_3 ] \rightarrow ( \ell_1^q , \ell_2^q , \ell_3^q )$ 
	\end{center}
	of the full plane $PG(2,q^2)$ to the biaffine plane of order $q^2$.    
	Let $G_{q^2}^{ \pi }$ be the corresponding polarity graph.  
	The vertex set of this graph is $\mathbb{F}_{q^2}^2$ and vertices
	$(p_1,p_2)$ and $(r_1,r_2)$ are adjacent if and only if 
	\begin{center}
		$p_2 + r_2^q = p_1 r_1^q$.
	\end{center}
	The graph $G_{q^2}^{ \pi}$ has $q^4 - q^3$ vertices of degree $q^2$,
	and the remaining $q^3$ vertices are absolute points and have degree $q^2 - 1$.  These absolute points are, up to isomorphism, contained in a classical unital in $PG(2,q^2)$, which has $q^3 +1$ points.    
	Our partition will have the property that each of these $q^3$ points
	will be in their own distinct part.

	\begin{proof}[Proof of Theorem \ref{thm:achromatic biaffine}] 
		Let $\beta \in \mathbb{F}_{q^2}$ be chosen so that $\{ \beta , \beta^q \}$ is a basis for $\mathbb{F}_{q^2}$ over $\mathbb{F}_q$.  It is a known fact in field theory that such a basis exists and it is typically called a \emph{normal basis}.      
		For $x \in \mathbb{F}_{q^2}$ and $y \in \mathbb{F}_q$, let 
		\[
		\mathcal{V}_{x,y}= \{ ( x , a \beta + y \beta^q ) : a \in \mathbb{F}_q \},
		\]
		and $\mathcal{V} = \{ \mathcal{V}_{x,y } : x \in \mathbb{F}_{q^2}, y \in \mathbb{F}_q \}$.  The family $\mathcal{V}$ forms a partition of the vertex set of $G_{q^2}^{ \pi } $.  Suppose   
		$\mathcal{V}_{x,y}$ and 
		$\mathcal{V}_{z,w}$ are two parts in $\mathcal{V}$, not necessarily distinct.  There will be an edge with one endpoint in the first part and the other in the second part if there exists $a,b \in \mathbb{F}_q$ such that 
		\begin{equation}\label{eq:c4 unitary 100}
			(a \beta + y \beta^q ) + ( b \beta + w \beta^q )^q
			=
			x z^q.  
		\end{equation}
		The left hand side can be rewritten as 
		\[
		a \beta + y \beta^q + b^q \beta^q + w^q \beta^{q^2} 
		= 
		a \beta + y \beta^q + b \beta^q + w^q \beta
		=
		(a + w)\beta + (y + b)\beta^q.
		\]
		Let $xz^q = s_{x,z} \beta + t_{x,z} \beta^q$ where $s_{x,z},t_{x,z} \in \mathbb{F}_q$.
		We then have 
		\[
		(a + w ) \beta + ( y + b ) \beta^q = s_{x,z} \beta + t_{ x,z} \beta^q
		\]
		which has the unique solution $a = s_{x,z} - w$ and
		$b= t_{x,z} - y$.  
		If the parts $\mathcal{V}_{x,y}$ and $\mathcal{V}_{z,w}$ are distinct, then 
		\[
		\{ ( x , ( s_{x,z}-w) \beta + y \beta^q ) , 
		(z , ( t_{x,z} - y ) \beta + w \beta^q ) \}
		\]
		is the unique edge between these parts.  
		Now assume  $\mathcal{V}_{x,y}=\mathcal{V}_{z,w}$, 
		so $x = z$ and $y = w$.  In this case, (\ref{eq:c4 unitary 100}) becomes
		\[
		(a \beta + y \beta^q ) + ( b \beta + y \beta^q)^q = x^{q+1}.
		\]
		Writing $x^{q+1} = s_x \beta + t_x \beta^q$ with $s_x, t_x \in \mathbb{F}_q$, we have
		\[
		s_x \beta + t_x \beta^q = x^{q+1} = (x^{q+1})^q = 
		(s_x \beta + t_x \beta^q )^q = t_x \beta + s_x \beta^q.
		\]
		This computation shows that $s_x = t_x$, and so we have 
		the unique solution $a = s_x -y  = b$.  This means that within the part $\mathcal{V}_{x,y}$ there is a loop at the vertex
		$( x , (s_x - y ) \beta + y \beta^q )$, but no edges.

		The partition $\mathcal{V}$ shows $\chi_a (G_{q^2}^{ \pi } ) \geq q^3$.  Observe that $G_{q^2}^{ \pi }$ has exactly $q^3$ vertices with loops.  Indeed, there are $q^3$ choices for $x$ and $q$ choices for $y$, and then the vertex $(x , ( s_x - y) \beta + y \beta^q )$
		is the only loop having $x$ in the first coordinate, and $y$ as the $\beta^q$ coefficient in the second coordinate.   The number of vertices of degree $q^2 - 1$ in $G_{q^2}^{ \pi }$ is $q^3$, and all other vertices have degree $q^2$.  Thus,
		\[
		e( G_{q^2}^{ \pi } ) = (1/2)[ (q^4 - q^3) q^2 + q^3( q^2-1) ] = \frac{q^6}{2} - \frac{q^3}{2} = \binom{q^3}{2}
		\]
		which shows this coloring is optimally complete, and  
		$ \chi_a (G_{q^2}^{ \pi } ) =q^3 =  \psi ( G_{q^2}^{ \pi } )$.
	\end{proof}
	
	\begin{proof}[Proof of Corollary \ref{cor:C4 partition theorem}]
		The graph $G_{q^2}^{ \pi }$ is a $C_4$-free graph.  The partition $\mathcal{V}$ from the proof of Theorem \ref{thm:achromatic biaffine} implies $f( q^3 , C_4) \leq q$ whenever $q$ is a prime power.  Applying a standard density of primes argument gives $f(r,C_4) \leq r^{1/3} +o(r^{1/3})$.  The lower bound $f(r,C_4) \geq (1 + o(1))r^{1/3}$ follows from the formula $\textup{ex}(N,C_4) = \frac{1}{2}N^{3/2} + o(N^{3/2})$ (see \cite{Axenovich} or \cite{Byrne} for further details).  
	\end{proof}

	
	\section{Proof of Theorem \ref{thm:c6 and c10}}\label{sec:C6 and C10}
	
	\subsection{A polarity graph of the biaffine generalized quadrangle $Q(4, q)$}
	
	Let $q$ be a power of a prime and let $GQ_q$ be the bipartite graph with parts $\mathcal{P}$ and $\mathcal{L}$, each of copy of 
	$\mathbb{F}_q^3$.  Vertex $(p_1,p_2,p_3) \in \mathcal{P}$ is adjacent to $[ \ell_1 , \ell_2 , \ell_3 ] \in \mathcal{L}$ if and only if 
	\begin{center}
		$p_2 + \ell_2 = p_1 \ell_1$ ~~ and ~~
		$p_3 + \ell_3 = p_1^2\ell_1$.
	\end{center}
	The graph $GQ_q$ is isomorphic to the incidence graph of 
	the biaffine part of the generalized quadrangle $Q(4, q)$.  
	It is a $q$-regular graph with girth 8 and $q^3$ vertices in each part.   
	The classical generalized quadrangle of order $q$ admits a polarity if $q = 2^{2e+1}$, $e \geq 1$ an integer \cite{Jacques}.  The same is true for the biaffine part.  
	
	Let $e \geq 1$ be an integer and $q = 2^{2e+1}$.  
	The function $\pi : V( GQ_q) \rightarrow V( GQ_q)$ defined by 
	\begin{center}
		$ \pi((p_1, p_2, p_3)) = [p_1^{2^{e+1}}, p_3^{2^e}, p_2^{2^{e+1}}] $ ~~and ~~
		$\pi([\ell_1, \ell_2, \ell_3]) = (\ell_1^{2^e}, \ell_3^{2^e}, \ell_2^{2^{e+1}})$
	\end{center}
	is a polarity.  Using the identity $x^{2^{2e+1}} = x$ for all $x \in \mathbb{F}_q$, it easy to check that $\pi^2$ is the identity map and by definition, $\pi$ interchanges $\mathcal{P}$ and $\mathcal{L}$.  For completeness, we verify that $\pi$ preserves adjacency in $GQ_q$.  The vertices $(x_1,x_2,x_3)$ and $[y_1,y_2,y_3]$ are adjacent if and only if 
	$x_2 + y_2 = x_1 y_1$ and $x_3 + y_3 = x_1^2 y_1$.  This pair of equations is equivalent to the pair 
	$(x_2 +y_2)^{2^{e+1}} = x_1^{2^{e+1}} y_1^{2^{e+1}}$ and
	$(x_3 + y_3)^{2^e} = (x_1^2)^{2^e} y_1^{2^e}$ which can be rewritten as
	\begin{center}
		$y_3^{2^e} + x_3^{2^e} = y_1^{2^e}  x_1^{2^{e+1}}$ and 
		$y_2^{2^{e+1}} + x_2^{2^{e+1}} = ( y_1^{2^{e} } )^2 x_1^{2^{e+1}}$. 
	\end{center}
	These hold if and only if $(y_1^{2^{e}} , y_3^{2^{e}} , y_2^{ 2^{e+1}} )$ is adjacent to 
	$[ x_1^{2^{e+1}} , x_3^{2^{e}} , x_2^{2^{e+1}}]$ in $GQ_q$.  
	We conclude that $\pi$ is a polarity and let 
	$GQ_q^{ \pi}$ be the corresponding polarity graph. 
	The vertex set of this graph is $\mathbb{F}_q^3$ and 
	distinct vertices $(p_1,p_2,p_3)$ and
	$(r_1,r_2,r_3)$ are adjacent if and only if  
	\begin{center}
		$p_2 + r_3^{2^e}= p_1r_1^{2^{e+1}}  $ ~~and ~~
		$  p_3 + r_2^{2^{e+1}} = p_1^2r_1^{2^{e+1}} $.  
	\end{center}
	The graph $GQ_q^{ \pi}$ has $q^3 - q^2$ vertices of degree $q$, and $q^2$ vertices of $q-1$, the latter of which are the absolute points of the polarity $\pi$.

	\begin{proof}[Proof of Theorem \ref{thm:c6 and c10} (a)]
		For $p_1,p_2 \in \mathbb{F}_q$, let $\mathcal{P}_{p_1, p_2} = \{(p_1, p_2, a): a \in \F_q \}$. The family 
		$\{ \mathcal{P}_{ p_1, p_2 } : p_1,p_2 \in \mathbb{F}_q \}$ partitions $V( GQ_q^{ \pi } )$ into 
		$q^2$ parts, each with $q$ vertices.  Suppose 
		that $\mathcal{P}_{p_1, p_2}$ and $\mathcal{P}_{r_1, r_2}$ are two parts in this partition, not necessarily distinct. 
		There is an edge between them if and only if there exists $a,b \in \mathbb{F}_q$ such that 
		\begin{center}
			$p_2 + b^{2^e} = p_1r_1^{2^{e+1}}$ ~~ 
			and ~~
			$a + r_2^{2^{e+1}} = p_1^2r_1^{2^{e+1}}$.
		\end{center}
		The first equation implies 
		$b = p_1^{2^{e+1}}r_1^2 +p_2^{2^{e+1}}$ while the second implies $a = p_1^2r_1^{2^{e+1}} + r_2^{2^{e+1}}$.
		This gives a unique solution for $a$ and $b$ which means 
		that there is exactly one edge, namely 
		\begin{equation}\label{eq:GQ edge between parts}
			\{ (p_1,p_2, p_1^2r_1^{2^{e+1}} + r_2^{2^{e+1}}) , 
			( r_1,r_2, p_1^{2^{e+1}}r_1^2 +p_2^{2^{e+1}} ) \},
		\end{equation} 
		with one endpoint in 
		$\mathcal{P}_{p_1,p_2}$ and $\mathcal{P}_{r_1,r_2}$. 
		If these parts are the same, so $p_1=r_1, p_2=r_2$, then 
		(\ref{eq:GQ edge between parts}) collapses to the single vertex $(p_1,p_2,p_1^{2^{e+1}+2} + p_2^{2^{e+1}})$.  
		This means that there is a loop at this vertex and no other edges in $\mathcal{P}_{p_1,p_2}$.  
		
		The partition $\{ \mathcal{P}_{p_1,p_2} : p_1,p_2 \in \mathbb{F}_q \}$ is an achromatic coloring with $q^2$ colors so $\chi_a ( GQ_q^{ \pi } ) \geq q^2$.  
		Also, any complete partition of $GQ_q^{ \pi }$ with more than $q^2$ parts must contain a part with at most $q-1$ vertices.  These $q-1$ vertices send out at most $(q-1) q = q^2 - q$ edges, but there are at least $q^2$ other parts so the partition is not complete.  We conclude that 
		$\chi_a ( GQ_q^{ \pi } ) = q^2= \psi ( GQ_q^{ \pi } ) $.  
	\end{proof}
	
	\noindent
	\textbf{Remark:} The proof that $\chi_a ( GQ_q^{ \pi } ) \leq q^2$
	is the same argument used in \cite{Ahlswede}, \cite{Maistrelli} for Paley graphs and hypercubes, respectively.  In general, if $G$ is an $n$-vertex graph with maximum degree $\Delta$ having a complete partition with $r$ parts, then $\lfloor n / r \rfloor \Delta \geq r - 1$.   Because this will be used again later, we state it as a proposition.
	
	\begin{prop}\label{prop}
		Let $G$ be an $n$-vertex graph with maximum degree $\Delta (G)$.  If $G$ has a complete partition with $r$ parts, then $\lfloor n / r \rfloor \Delta (G) \geq r - 1$. 
	\end{prop}

	
	\subsection{A polarity graph of the biaffine split Cayley hexagon}
	
	Let $q$ be a power of a prime and let $\Gamma_q$ be the bipartite graph with parts $\mathcal{P}$ and $\mathcal{L}$, each a copy of $\mathbb{F}_q^5$.  Vertex
	$(p_1, p_2 , p_3,p_4,p_5) \in \mathcal{P}$ is adjacent to vertex $[ \ell_1 , \ell_2 , \ell_3 , \ell_4 , \ell_5 ] \in \mathcal{L}$ if the following equations are satisfied: 
	\begin{align*}
		p_2 + \ell_2 &= p_1 \ell_1, \\
		p_3 + \ell_3 &= p_1\ell_2, \\
		p_4 + \ell_4 &= p_1\ell_3, \\
		p_5 + \ell_5 &= p_2\ell_3 - p_3\ell_2.
	\end{align*}
	The graph $\Gamma_q$ is isomorphic to the incidence 
	graph of the biaffine part of the split Cayley hexagon (\cite{LSW} attributes these equations to Jason Williford).  It is a $q$-regular graph with girth 12.  
	When $q  = 3^{ \eta }$, these equations can be re-written in terms of just $p_1$ and $\ell_1$, as observed by Williford \cite{Williford}.
	
	\begin{lem}
		Let $e\geq 1$ and $q = 3^e$. Then the graph $\Gamma_q$ defined above is isomorphic to the algebraically defined graph $GH_q$ constructed over $\F_q^5$ with equations
		\begin{equation}\label{eq:GHQ3 equations}
			p_2 + \ell_2 = p_1 \ell_1, ~~~
			p_3 + \ell_3 = p_1^2\ell_1, ~~~
			p_4 + \ell_4 = p_1^3\ell_1,~~~
			p_5 + \ell_5 = p_1^3\ell_1^2.
		\end{equation}
	\end{lem}
	
	\begin{proof}
		Define $\phi: V( \Gamma_q ) \rightarrow V( GH_q)$ by 
		\begin{align*}
			\phi((p_1, p_2, p_3, p_4, p_5)) &= (p_1, p_2, p_3 + p_1p_2, p_4 + p_1p_3 + p_1^2p_2, -p_5 + p_2^2p_1 - p_2p_3) = (p_1', p_2', p_3', p_4', p_5') \\
			\\
			\phi([\ell_1, \ell_2, \ell_3, \ell_4, \ell_5]) &= [\ell_1, \ell_2, \ell_3, \ell_4, -\ell_5 + \ell_2\ell_3] = [\ell_1', \ell_2', \ell_3', \ell_4', \ell_5'].
		\end{align*}
		We will show that $\phi$ is an isomorphism from $\Gamma_q$ to $GH_q$.
		It is clear that $\phi$ is a bijection. We have 
		\begin{align*}
			p_2' + \ell_2' &= p_2 + \ell_2 = p_1\ell_1 = p_1'\ell_1', \\
			p_3' +\ell_3' &= p_3 +p_1p_2 + \ell_3 = p_1\ell_2 + p_1p_2 = p_1^2\ell_1 = (p_1')^2\ell_1' , \\
			p_4' + \ell_4' &= p_4 + p_1p_3 + p_1^2p_2 + \ell_4 =  p_1p_3 + p_1^2p_2 + p_1\ell_3 = p_1(p_3 + \ell_3) +p_1^2p_2 = p_1^2\ell_2+ p_1^2p_2 = (p_1')^3\ell_1', \\
			p_5' + \ell_5' &= -p_5 + p_2^2p_1 - p_2p_3 -\ell_5 + \ell_2\ell_3 = 2p_2\ell_3 + p_3\ell_2 +p_2^2p_1 +2p_2p_3 + \ell_2\ell_3 \\
			&= (p_2+\ell_2)(p_3+\ell_3) +p_2\ell_3 + p_2^2p_1+p_2p_3 = (p_2+\ell_2)(p_3+\ell_3) +p_2(p_1\ell_2 - p_3)+ p_2^2p_1+p_2p_3 \\
			&= (p_2+\ell_2)(p_3 + p_1p_2 + \ell_3) = p_1^3\ell_1^2 = (p_1')^3(\ell_1')^2.
		\end{align*}
		Thus, the adjacency relations in $GH_q$ hold and so $\phi$ is an isomorphism.
	\end{proof}

	From here on, assume $q = 3^{2e + 1}$.  The function $\pi : V (GH_q) \rightarrow V( GH_q)$ defined by 
	\begin{eqnarray*}
		\pi((p_1, p_2, p_3, p_4, p_5)) &=& [p_1^{3^{e+1}}, p_4^{3^e} ,p_5^{3^e}, p_2^{3^{e+1}}, p_3^{3^{e+1}}] ~\mbox{and} \\
		\pi([\ell_1, \ell_2, \ell_3, \ell_4, \ell_5]) &= &(\ell_1^{3^{e}}, \ell_4^{3^e} ,\ell_5^{3^e}, \ell_2^{3^{e+1}}, \ell_3^{3^{e+1}})
	\end{eqnarray*}
	is a polarity of $GH_q$.  As in the case of $GQ_q$, it is straightforward to verify that $\pi$ interchanges the parts and $\pi^2$ is the identity map.  The verification that 
	$\pi$ preserves adjacency is similar to $GQ_q$.

	Let $GH_q^{ \pi }$ be the corresponding polarity graph.  The vertex set is $\mathbb{F}_q^5$ where $q = 3^{2e + 1}$.  Vertices $(p_1, p_2,p_3,p_4 , p_5)$ and 
	$(r_1 , r_2,r_3,r_4 , r_5 )$ are adjacent if and only if 
	\begin{equation}\label{eq:GHQ3 polarity equations}
		p_2 + r_4^{3^e} = p_1 r_1^{3^{e+1}}, ~~~
		p_3 + r_5^{3^e} = p_1^2 r_1^{3^{e+1}} , ~~~
		p_4 + r_2^{3^{e+1} } = p_1^3 r_1^{e^{e+1} } ,~~~
		p_5 + r_3^{3^{e+1} } = p_1^3 r_1^{ 2 \cdot 3^{e+1} } .
	\end{equation}
	
	\begin{proof}[Proof of Theorem \ref{thm:c6 and c10} (b)]
		For $p_1,p_2,p_3 \in \mathbb{F}_q$, let 
		\[
		\mathcal{P}_{p_1,p_2,p_3} = \{ ( p_1 , p_2 , p_3 , a , b ) : a,b \in \mathbb{F}_q \}.
		\]
		The family $\{ \mathcal{P}_{ p_1 ,p_2,p_3 } : p_1,p_2,p_3 \in \mathbb{F}_q \}$ partitions $V( GH_q^{ \pi } )$ into $q^3$ parts, each with $q^2$ vertices.  There is an edge between 
		$\mathcal{P}_{p_1,p_2,p_3}$ and 
		$\mathcal{P}_{r_1,r_2,r_3}$ if and only if there exists $a,b,c,d \in \mathbb{F}_q$ such that 
		\begin{center}
			$p_2 + c^{3^e} = p_1 r_1^{3^{e+1} } $, ~
			$p_3 + d^{3^e} = p_1^2 r_1^{3^{e+1} } $, ~
			$a + r_2^{3^{e+1}} = p_1^3 r_1^{3^{e+1}}$, ~
			and ~
			$b + r_3^{3^{e+1}} = p_1^3 r_1^{ 2 \cdot 3^{e+1} }$.
		\end{center}
		The unique solution for $a,b,c,d$ is 
		\begin{center}
			$a = p_1^3 r_1^{3^{e+1} } - r_2^{3^{e+1}} $, 
			$b = p_1^3 r_1^{2 \cdot 3^{e+1} } - r_3^{3^{e+1}}$, 
			$c = ( p_1 r_1^{3^{e+1} } - p_2 )^{3^{e+1} } $, 
			and $d = ( p_1^2 r_1^{3^{e+1} } - p_3)^{3^{e+1} }$.    
		\end{center}
		This means there is exactly one edge between $\mathcal{P}_{p_1,p_2,p_3}$ and 
		$\mathcal{P}_{r_1,r_2,r_3}$.  If these two parts are the same, then the unique solution has 
		\[
		c = ( p_1 p_1^{3^{e+1}} - p_2)^{3^{e+1}} 
		=
		p_1^{3^{e+1}} p_1^3 - p_2^{3^{e+1}} = 
		p_1^3 p_1^{3^{e+1} } - p_2^{ 3^{e+1}} = a
		\]
		and
		\[
		d = ( p_1^2 p_1^{3^{e+1}} - p_3 )^{ 3^{e+1}} 
		= p_1^{2 \cdot 3^{e+1}} p_1^3 - p_3^{3^{e+1}} = b.
		\]
		Thus, within any part $\mathcal{P}_{p_1,p_2,p_3}$, there is no edge and only a loop at the vertex 
		\[
		(p_1,p_2,p_3 , p_1^3 p_1^{3^{e+1}} - p_2^{3^{e+1}} , p_1^3 p_1^{2 \cdot 3^{e+1}} - p_3^{3^{e+1}}  ).
		\]
		The partition $\{ \mathcal{P}_{p_1,p_2,p_3} : p_1,p_2,p_3 \in \mathbb{F}_q \}$ shows $\chi_a ( GH_q^{ \pi } ) \geq q^3$.  By Proposition \ref{prop}, there is no complete partition into $q^3+1$ parts since $\lfloor q^5 / (q^3+1) \rfloor q < (q^3-1)+1$.     We conclude that $\chi_a ( GH_q^{ \pi } ) =q^3 =  \psi ( GH_q^{ \pi} )$.  
	\end{proof}

	\begin{proof}[Proof of Corollary \ref{cor:C6 and C10}] 
		Lam and Verstra\"{e}te \cite{Lam} proved that for $m \geq 2$, 
		one has 
		\[
		\textup{ex}(n, , \{C_4, C_6 , \dots , C_{2 m} \} ) \le \frac{1}{2} n^{1 + 1/ m } + 2^{ m^2}  
		~~\mbox{for all $n$.}
		\]
		They noted that the polarity graphs of generalized $(m+1)$-gons have $\frac{1}{2}m^{1+1/k} + O(m)$ edges given a matching lower bound up to an error term.  
		This upper bound and (\ref{ineq:trivial turan ub}) implies 
		\[
		f( r , \{ C_4 , C_6 \} ) \geq (1 + o(1)) r^{1/2}.
		\]
		On the other hand, for any $r$ of the form $r = 2^{2(2e+1)} $ where $e \geq 1$ is an integer, we have $f( r , \{C_4 , C_6 \} ) \leq 2^{2e+1}$ because
		$GQ_q^{ \pi }$ is a $\{C_4 ,C_6 \}$-free $(r,r^{1/2})$-graph.    
		
		A similar argument applies to the $\{C_4,C_6,C_8 , C_{10} \}$-free 
		$(q^3,q^2)$-graph $GH_q^{ \pi}$, $q= 3^{2e + 1}$.  
	\end{proof}

	\section{Optimally complete polarity graphs and general pseudocolorings}

	Wenger \cite{Wenger} discovered an important family of algebraically defined graphs.  Using these graphs, it was proved  in \cite{Byrne} that $f(r,C_6) \leq 2r^{1/2} + o(r^{1/2})$ and $f(r,C_{10}) \leq 2r^{2/3} + o(r^{2/3} )$. 
	This is done by defining complete partitions of the $C_6$-free and $C_{10}$-free graphs from \cite{Wenger}.  
	With some care, this approach together with the ideas of Section \ref{sec:C4} can be used to prove a lower bound that applies to all algebraically defined graphs.  First we need a definition. Let $m \geq 2$ be an integer and $q$ be a power of a prime.  Functions $f_2 , f_3, f_4, \dots , f_m$ with $f_i : \mathbb{F}_{q^2}^{2i-2} \rightarrow \mathbb{F}_{q^2}$ 
	are called \emph{point-line-symmetric} if  
	\[
	f_j ( x_1 , y_1 , x_2 , y_2 , \dots , x_{j-1} , y_{j-1} ) 
	= 
	f_j ( y_1 , x_1 , y_2 , x_2 , \dots , 
	y_{j-1} , x_{j-1} ) 
	\]
	for all $x_1,x_2, \dots , x_{j-1} ,y_1,y_2, \dots , y_{j-1} \in \mathbb{F}_{q^2}$ and $2 \leq j \leq m$.
	
	\begin{thm}\label{theorem:general upper bound}
		Let $m \geq 2$ be an integer and $q$ be a power of a prime.
		Let $\rho = q^2$ if $m$ is even and $\rho = q$ if $m$ is odd.  If $\Gamma_{ \rho }$ is any algebraically defined graph over $\mathbb{F}_{ \rho }^m$, then 
		\[
		\psi ( \Gamma_{ \rho } ) \geq \sqrt{ e ( \Gamma_{ \rho } ) }.
		\]
	If $m$ is even and the functions (polynomials) used to defined adjacency are point-line-symmetric with coefficients in $\F_q$, then there is a polarity $\pi$ of $\Gamma_{q^2}$ such that 
	$\chi_a ( \Gamma_{q^2}^{ \pi } ) = \sqrt{ 2 e ( \Gamma_{q^2}^{ \pi } ) +1/4 } +1/2 $, i.e., $\Gamma_{q^2}^{ \pi}$ is optimally complete.  
	\end{thm}
	\begin{proof} 
		First assume that $m \geq 3$ is odd.  Let $q$ be a power of a prime and $\Gamma_q$ be the algebraically defined graph over $\F_q^{m}$ with  adjacency equations
		\begin{center}
			$\ell_j + p_j = f_j ( \ell_1 , p_1 , \ell_2 , p_2 , \dots , \ell_{j-1} , p_{j-1} )$, ~~$j = 2,3, \dots ,m $.
		\end{center}
		For $p_1,p_3,p_5,\dots,p_m \in \mathbb{F}_q$ and 
		$\ell_1,\ell_2,\ell_4,\dots , \ell_{m-1}  \in \mathbb{F}_q$, let 
		\begin{align*}
			\mathcal{P}_{p_1,p_3,p_5,\dots , p_m } 
			&= \{ (p_1 , p_2 , \dots , p_m) \in \mathcal{P} : p_2, p_4 , \dots , p_{m-1} \in \mathbb{F}_{q} \} , \\
			\mathcal{L}_{ \ell_1 , \ell_2 , \ell_4 , \dots , \ell_{m-1} } 
			&= 
			\{ [ \ell_1 , \ell_2 , \dots , \ell_m ] \in \mathcal{L} : \ell_3 ,\ell_5 , \dots , \ell_m \in \mathbb{F}_q \}. 
		\end{align*}
		Given any $\mathcal{P}_{p_1,p_3,p_5,\dots , p_m } $ 
		and $\mathcal{L}_{ \ell_1 , \ell_2 , \ell_4 , \dots , \ell_{m-1} } $, there will be exactly one edge between these two sets. 
		Indeed, suppose that $p_1,p_3,p_5 , \dots , p_m$ and 
		$\ell_1 , \ell_2 , \ell_4  , \dots , \ell_{m-1}$ are all fixed elements of $\mathbb{F}_q$.  
		The $j=2$ equation $\ell_2 + p_2 = f_2( \ell_1  , p _1)$ has a unique solution for $p_2$ since 
		$p_1,\ell_1,\ell_2$ are fixed.  
		Once $p_2$ is determined, the $j=3$ equation
		$\ell_3 +p_3 = f_3 ( \ell_1 , p_1 , \ell_2 , p_2 )$
		will have a unique solution for $\ell_3$ since $p_2,p_3,p_1 , \ell_2,\ell_1$ are all determined.  Likewise, if $m > 3$, then the $j = 4$ equation $\ell_4 + p_4 = f_4 ( \ell_1 , p_1 , \ell_2 , p_3 , \ell_3 , p_3)$ will have a unique solution for $p_4$ since $p_3,p_2,p_1,\ell_3,\ell_2,\ell_1$ are now all determined.  This process continues and at end, the $j = m$ equation will give $\ell_m$.  
		There are $q^{ \frac{m-1}{2} } $ possible sets
		$\mathcal{P}_{p_1,p_3,p_5,\dots , p_m } $ and the same is true for the number of $\mathcal{L}_{ \ell_1 , \ell_2 , \ell_4 , \dots , \ell_{m-1} } $'s. 
		Arbitrarily form a pairing between the
		$\mathcal{P}_{p_1,p_3,p_5,\dots , p_m } $'s  
		and $\mathcal{L}_{ \ell_1 , \ell_2 , \ell_4 , \dots , \ell_{m-1} } $'s.  The union of each of the pairs 
		form the parts of a complete partition of $\Gamma_q$ into
		$q^{ \frac{m+1}{2} }$ parts.  Therefore, 
		$\psi ( \Gamma_q ) \geq q^{ \frac{m+1}{2} }$.  
		
		Now assume that $m \geq 2$ is even.  Let $\Gamma_{q^2}$ be the algebraically defined graph with the same functions used to define adjacency, except now over $\mathbb{F}_{q^2}$.  Let $\{ 1 , \mu \}$ be a $\mathbb{F}_q$-basis for $\mathbb{F}_{q^2}$. For $p_1,p_3,p_5 , \dots , p_{m-1} \in \mathbb{F}_{q^2}$, $p_{m}' \in \mathbb{F}_q$
		and 
		$\ell_1,\ell_2,\ell_4 , \dots , \ell_{m-2} \in \mathbb{F}_{q^2}$, $\ell_{m}'' \in \mathbb{F}_q$, let
		\begin{align*}
			\mathcal{P}_{p_1,p_3,p_5,\dots , p_{m-1} , p_{m}' } 
			&= \{ (p_1 , p_2 , \dots , p_{m-1} , p_m' + p_m'' \mu ) \in \mathcal{P} : p_2, p_4 , \dots , p_{m-2} \in \mathbb{F}_{q^2} , p_{m}'' \in \mathbb{F}_q \}, \\
			\mathcal{L}_{ \ell_1 , \ell_2 , \ell_4 , \dots , \ell_{m-2}, \ell_{m}'' } 
			&= 
			\{ [ \ell_1 , \ell_2 , \dots , \ell_{m-1},  \ell_m ' + \ell_m '' \mu  ] \in \mathcal{L}  : \ell_3 ,\ell_5 , \dots , \ell_{m-1}  \in \mathbb{F}_{q^2} , \ell_{m}' \in \mathbb{F}_q \}. 
		\end{align*}
		Let us now show that between any two of these sets, there is exactly one edge.  Fix a 
		$\mathcal{P}_{p_1,p_3,p_5,\dots , p_{m-1} , p_{m}' }$
		and an 
		$\mathcal{L}_{ \ell_1 , \ell_2 , \ell_4 , \dots , \ell_{m-2}, \ell_{m}'' } $.  As in the case when $m$ is odd, the $j = 2$ equation will determine $p_2$, the $j = 3$ equation will determine $\ell_3$, and so on.  The $j = m-1$ equation will determine $\ell_{m-1}$ and it remains to examine the $j  = m$ equation  
		\[
		( \ell_m' + \ell_{m}'' \mu ) + 
		( p_m' + p_{m}'' \mu ) =
		f_m (  \ell_1 , p_1 , \ell_2 , p_2 , \dots , \ell_{m-1} , p_{m-1} ).  
		\]
		This can be rewritten as 
		\[
		( \ell_m' + p_m') + ( \ell_m'' + p_m'' ) \mu 
		= f_m (  \ell_1 , p_1 , \ell_2 , p_2 , \dots , \ell_{m-1} , p_{m-1} ).
		\]
		Because $\ell_m''$, $ p_m'$, and the right hand side are all fixed, 
		we have a unique solution for $\ell_m'$ and $p_m''$.  
		The remaining steps of the argument are the same as the case when $m$ is odd.    
		The result is the lower bound $\psi ( \Gamma_{q^2} ) \geq q^{m+1}$.  
		
		Continuing with assumption that $m$ is even, suppose further that the functions $f_2, f_3, \dots , f_m$ defining adjacency are point-line symmetric.  
		\[
		f_j (x_1, y_1 , x_2, y_2, \dots ,x_{j-1}, y_{j-1})
		= f_j (y_1, x_1 , y_2, x_2 , \dots , y_{j-1} , x_{j-1} ). 
		\]
		Define $\pi : V( \Gamma_{q^2} ) \rightarrow V( \Gamma_{q^2} )$ by the rule
		\[
		\pi ( ( p_1 ,p_2 , \dots , p_m )) 
		= 
		[ p_1^q , p_2^q , \dots , p_m^q]
		~~\mbox{and}~~
		\pi ([ \ell_1 , \ell_2 , \dots , \ell_m ] ) 
		= 
		( \ell_1^q , \ell_2^q , \dots , \ell_m^q ).  
		\]
		Certainly $\pi$ is a bijection that interchanges $\mathcal{P}$ and $\mathcal{L}$.  Suppose
		that $(p) = (p_1, \dots , p_m)$ is adjacent to 
		$(r) = ( r_1 , \dots , r_m )$ in $\Gamma_{q^2}$. By definition, this implies that
		\[
		r_j + p_j = f_j ( r_1 , p_1 , \dots , r_{j-1}, p_{j-1}) ~~\mbox{for $j=2,3,\dots , m$}.  
		\]
		Observe that $\pi((p))$ must be adjacent to $\pi((r))$ since for each $j$, 
		\begin{eqnarray*} 
			p_j^q + r_j^q = ( r_j  + p_j)^q &=& 
			f_j ( r_1 , p_1 , \dots , r_{j-1} , p_{j-1})^q \\
			& = & 
			f_j ( p_1 , r_1, \dots , p_{j-1} , r_{j-1} )^q \\
			& = & f_j ( p_1^q , r_1^q, \dots , p_{j-1}^q , r_{j-1}^q ).
		\end{eqnarray*}
		This computation shows that $\pi((r_1 , \dots , r_m)$ is 
		adjacent to $\pi (( p_1 , \dots , p_m ))$.  
		We conclude that $\pi$ is a polarity.  Let 
		$\Gamma_{q^2}^{ \pi}$ be the polarity graph of $\Gamma_{q^2}$ with respect to $\pi$. The vertex set of $\Gamma_{q^2}^\pi$ is $\mathbb{F}_{q^2}^m$.  Vertices $(p_1,p_2, \dots , p_m)$ and $(r_1,r_2, \dots , r_m)$ are adjacent if 
		\[
		r_i^q + p_i =f_i ( r_1^q , p_1 , r_2^q , p_2 , \dots , r_{i-1}^q , p_{i-1} )  \mbox{~~for $i=2,3, \dots , m$}.
		\]
		For $x_1 \in \mathbb{F}_{q^2}$ and $y_2,y_3, \dots , y_{m} \in \mathbb{F}_q$, define
		\[
		\mathcal{V}_{x_1, y_2, y_3, \dots , y_{m} } 
		= 
		\{ 
		( x_1 , a_2 \beta + y_2 \beta^q , 
		a_3 \beta + y_3 \beta^q , \dots , 
		a_{m} \beta + y_{m} \beta^q ) :
		a_2,a_3, \dots ,a_{m} \in \mathbb{F}_q \}. 
		\]
		These sets partition $V(\Gamma_{q^2}^{ \pi})$ into $q^{m+1}$ parts.  Fix $x_1,z_1 \in \mathbb{F}_{q^2}$ and $y_2, \dots , y_{m} , w_2 , \dots , w_{m} \in \mathbb{F}_q$, and 
		consider $\mathcal{V}_{x_1, y_2 , y_3 , \dots ,y_{m} }$ and $\mathcal{V}_{ z_1, w_2 , w_3 , \dots , w_{m} }$.  
		Assume first these parts are distinct.  The 
		$i=2$ equation is 
		$ ( b_2 \beta + w_2 \beta^q)^q + ( a_2 \beta+ y_2 \beta^q ) = f_2 (z_1^q,x_1)$
		which is equivalent to 
		\[
		(a_2 + w_2 ) \beta + ( b_2 + y_2 ) \beta^q = f_2 (z_1^q , x_1).
		\]
		Since $x_1,z_1,y_2,w_2$ are all fixed, this equation has a unique solution for $a_2$ and $b_2$.  The $i=3$ equation is equivalent to 
		\[
		(a_3 + w_3) \beta + ( b_3 + y_3) \beta^q = f_3 (z_1^q , x_1 , ( b_2 \beta + w_2 \beta^q)^q , a_2 \beta + y_2 \beta^q ).
		\]
		This has a unique solution for $a_3$ and $b_3$ when all other values are fixed.  Like in the previous cases, this process continues and the result is a unique solution for $a_2, \dots , a_{m} , b_2 , \dots b_{m}$ which means that there is exactly one edge between these two parts.
		
		Now assume these parts are the same so that $x_1=z_1$ and $y_i = w_i$ for $2 \leq i \leq m$.  In this case, the $i=2$ equation is $(a_2 + y_2) \beta  +( b_2 + y_2) \beta^q = f_2 (x_1^q,x_1)$.  
		Using the point-line-symmetry of $f_2$, 
		\begin{eqnarray*}
			(b_2 + y_2) \beta + ( a_2 + y_2) \beta^q & = &  
			[  (a_2 +y_2) \beta + (b_2 +y_2) \beta^q ]^q \\
			& = &f_2 (x_1^q , x_1 )^q  = f_2 (x_1^q , x_1) \\
			& = &  (a_2 + y_2) \beta  +( b_2 + y_2) \beta^q.
		\end{eqnarray*} 
		Therefore, $b_2 +y_2 = a_2 + y_2$ so $b_2 = a_2$.  
		For $3 \leq i \leq m$, the $i$-th equation can be written as 
		\[
		(a_i + y_i) \beta + ( b_i + y_i ) \beta^q = 
		f_i ( x_1^q , x_1 , (a_2 \beta + y_2 \beta^q)^q , 
		a_2 \beta + y_2 \beta^q , 
		\dots 
		, 
		(a_{i-1} \beta + y_{i-1} \beta^q)^q , 
		a_{i-1} \beta + y_{i-1} \beta^q ).
		\]
		Again using the point-line-symmetry of $f_i$, we get $(b_i+y_i) \beta + ( a_i + y_i) \beta^q = (a_i + y_i) \beta + ( b_i +y_i) \beta^q$ so $b_i = a_i$.  The only edge in   
		$\mathcal{V}_{x_1, y_2, y_3, \dots , y_{m} } $
		is a loop at the vertex $(x_1,a_2,a_3,\dots ,a_m)$ where the $a_i$'s are uniquely determined by $x_1,y_2,\dots,y_m$ through the adjacency equations.  Since this partition has no edge within a part and exactly one edge between any two parts, it is an achromatic coloring and even more, $\Gamma_{q^2}^{ \pi}$ is optimally complete.    
	\end{proof} 
	
	Observe that since $\Gamma_q$ has $q^{m+1}$ edges and $\Gamma_{q^2}$ has $q^{2m+2}$ edges, we have
	\begin{equation}\label{eq:conclusion}
		\frac{\psi(\Gamma_q)}{\sqrt{2e(\Gamma_q)}} \geq \frac{1}{\sqrt{2}} \mbox{~~and~~} 
		\frac{ \psi ( \Gamma_{q^2} ) }{ \sqrt{ 2 e ( \Gamma_{q^2} ) } } \geq \frac{1}{ \sqrt{2}}.
	\end{equation}

	\section{Concluding Remarks}\label{section:conclusion} 
	
	Let us conclude by returning to the 
	problem of Ahlswede et al.\ \cite{Ahlswede} of finding explicit sequences with $\underline{\psi}(G) = 0$, but taking into account some of the 
	developments since then.  
	
	Cairnie and Edwards \cite{Cairnie1} proved that if 
	$\mathcal{G} = (G_1,G_2, \dots )$ is a sequence where $\Delta (G_i ) \leq C$ for some constant $C$, then $\mathcal{G}$ is almost optimally complete.  
	For such a sequence, we have $\underline{\psi} ( \mathcal{G} ) = \underline{\chi_a} ( \mathcal{G} ) = 1$.  So, as also noted by Roichman \cite{Roichman}, it is reasonable to assume the degrees of the graphs are unbounded.      
	One may also need an assumption on how close the graphs are to being regular.  For example, if $1 \leq n_1 < n_2$, then $\psi ( K_{n_1 ,n_2 } ) = n_1 + 1$.  Hence, the sequence $\mathcal{G} = \{ K_{i,i^2} \}_{n=1}^{ \infty }$ will satisfy $\underline{ \psi } ( \mathcal{G} ) = 0$.  One can construct similar sequences using complete $r$-partite graphs where the largest part size becomes larger and larger than the union of all the other parts.  The formulas for the psuedochromatic numbers of complete $r$-partite graphs can be found in \cite{Yeg3}.  
	
	Next, suppose that $H$ is a bipartite graph that contains a cycle.  In cases in which a good lower bound on  $\textup{ex}(n,H)$ is known, the constructions are often algebraically defined graphs.  As algebraically defined graphs have large pseudochromatic number, one could ask if a similar property holds for $H$-free extremal graphs.  This naturally suggests a connection to the research on $f(r,H)$ initiated by Axenovich and Martin \cite{Axenovich}. 
	They conjectured that there is a positive constant $c$, depending on $H$, such that if there is an $H$-free graph with $rk$-vertices and $cr^2$ edges, then $f(r,H) \leq k$.  Now, $f(r,H) \leq k$ is true if and only if there is an $H$-free $(r,k)$-graph $G$.  For this graph $G$, we have $\psi (G) \geq r$ and because the color classes do not have to be independent sets, we may assume $e(G)  =c r^2$.  Hence, $\frac{ \psi (G) }{  \sqrt{ 2 e(G)  } } \geq \frac{1}{ \sqrt{2c} }$.  The conclusion is that if their conjecture is true, then $H$-free extremal graphs would also not be the best place to search for examples that solve the 
	problem of Ahlswede et al.  Their conjecture has been confirmed in several special cases such as $C_4$, $C_6$, $C_{10}$, and $K_{t,(t-1)!+1}$ \cite{Axenovich, Byrne}.  
	Of course, by Theorem \ref{theorem:general upper bound} and the inequalities (\ref{eq:conclusion}), any sequence $\mathcal{G}$ of algebraically defined graphs or polarity graphs will satisfy $\underline{ \psi } (\mathcal{G} ) \geq 1/\sqrt{2} >  0$.  
	
	Finally, Roichman \cite{Roichman} posed the problem of finding sequences of Cayley graphs with large achromatic number, i.e., $\underline{\chi_a} ( \mathcal{G} ) > c $ for some positive constant $c$.  Since $\chi_a (G) \leq \psi (G)$, this would also give graphs with large pseudochromatic number.  
	There are many Cayley graphs that are not algebraically defined graphs, and conversely, not every algebraically defined graph is a Cayley graph.  
	A possible direction for further research is a combination of these ideas where one looks to study explicit Cayley graphs with small pseudochromatic number, or perhaps more generally, graphs satisfying a regularity condition on the vertex degrees.



\begin{thebibliography}{99}
		
		\bibitem{Ahlswede}
		R.\ Ahlswede, S.\ L.\ Bezrukov, 
		A.\ Blokhuis, K.\ Metsch, G.\ E.\ Moorhouse, 
		Partitioning the $n$-cube into sets with mutual distance 1,
		{\em Appl.\ Math.\ Lett}. Vol.\ 6, No.\ 4 (1993), 17--19.  
		
		\bibitem{Axenovich}
		M.\ Axenovich, R.\ Martin,
		Splits with forbidden subgraphs, 
		{\em Discrete Math}. 345(2), (2022), Article 112689.  
		
	
		\bibitem{Barbanera}
		F.\ Barbanera, 
		Covering numbers in Ramsey problems, 
		Master's thesis, Karlsruhe Institute of Techology, 2019. 
		
		\bibitem{Bollobas}
		B.\ Bollob\'{a}s, P.\ Catlin, P.\ Erd\H{o}s, 
		Hadwiger's conjecture is true for almost every graph,
		{\em Europ.\ J.\ Combin.\ Theory}, Vol.\ 1,  Issue 3, (1980), 195--199.  
		
		\bibitem{Byrne}
		J.\ Byrne, M.\ Tait, C.\ Timmons,
		Forbidden subgraphs and complete partitions,
		arXiv:2308:16728 2023.
		
		\bibitem{Cairnie1}
		N.\ Cairnie, K.\ Edwards, 
		Some results on the achromatic number, 
		{\em J.\ Graph Theory}, 26, No.\ 3 (1997), 129--136. 
		
		\bibitem{EdwardsTrees}
		K.\ J.\ Edwards,
		Achromatic number versus pseudoachromatic number:
		a counterexample to a conjecture of Hedetniemi,
		{\em Discrete Math}.\ 219 (2000), 271--274.  
		
		\bibitem{Kostochka}
		A.\ V.\ Kostochka,
		On the minimum of the Hadwiger number for graphs with a 
		given mean degree of vertices,
		{\em Amer.\ Math.\ Soc.\ Transl.} (2) Vol.\ 132, 1986, 15--31.  
		
		\bibitem{Lam}
		T.\ Lam, J.\ Verstra\"{e}te, 
		A note on graphs without short even cycles, 
		{\em Electron.\ J.\ Combin}.\ Vol 12 (2005), N5.  
		
		\bibitem{LSW}
		F. Lazebnik, S. Sun, Y. Wang, {\em Some families of graphs, hypergraphs and digraphs
			defined by systems of equations: a survey}, Selected topics in graph theory and
		its applications: Lecture Notes in Seminal Interdisciplinary Mathematics, 14, (2017): 105–142.
		
		\bibitem{LUW}
		F. Lazebnik, V.A. Ustimenko, A. Woldar, Polarities and $2k$-cycle free graphs, {\em Discrete Math.} 197-8, (1999), 503--513.
		
		\bibitem{Yeg3} 
		B.\ Logeshwary, V.\ Yegnanarayanan,
				On Pseudocomplete and Complete Coloring of Graphs. {\em Proceedings of the International Conference on Applied Mathematics and Theoretical Computer Science}, 2013.
		
		
		\bibitem{Maistrelli}
		E.\ Maistrelli, D.\ B.\ Penman,
		Some colouring problems for Paley graphs,
		{\em Discrete Math}. 306 (2006), 99--106.  
		
			
		\bibitem{Roichman}
		Y.\ Roichman,
		On the achromatic number of hypercubes,
		{\em J.\ Combin.\ Theory Ser.\ B} {\bf 79} (2000), 177--182.  
		
		\bibitem{Jacques}
		J.\ Tits, Ovo\"{i}des et groupes de Suzuki,
		{\em  Arch. Math.} 13 (1962), 187--198.
		
		\bibitem{Thomason}
		A.\ Thomason,
		Complete minors in pseudorandom graphs,
		{\em Random Structures Algorithms} 17 (1) (2000), 26--28.  
		
		\bibitem{Wenger}
		R.\ Wenger, Extremal graphs with no $C^4$'s, $C^6$'s, or $C^{10}$'s,
		{\em J.\ Combin.\ Theory Ser.\ B}, 52 (1991), 113--116.  
		
		\bibitem{Williford}
		J.\ Williford, \textit{personal communication}, (2022).
		
		\bibitem{Yeg1}
		V.\ Yegnanarayanan, 
		The pseudoachromatic number of a graph,
		{\em Southeast Asian Bull.\ Math}. 24 (2000), 129--136.  
		
		\bibitem{Yeg2}
		V.\ Yegnanarayanan,
		Graph colourings and partitions, 
		{\em Theoret.\ Comput.\ Sci}. 263 (2001), 59--74.  
	
		
		
	\end{thebibliography}
\end{document}